\documentclass{amsart}

\usepackage{amsmath}
\usepackage{amsthm}
\usepackage{amssymb}

\newcommand{\set}[2]{\left\{ #1 \mid #2 \right\}}
\renewcommand{\setminus}{-}
\newcommand{\numbers}{\mathbb{N}}
\newcommand{\integers}{\mathbb{Z}}

\newcommand{\sgn}[1]{\operatorname{sgn}(#1)}
\newcommand{\mdeg}[1]{\operatorname{deg}(#1)}

\newcommand{\comps}{c}
\newcommand{\components}[1]{\comps(#1)}

\newcommand{\djoin}{\cdot}

\newcommand{\lcm}{\operatorname{lcm}}

\newcommand{\boldx}{\mathbf{x}}
\newcommand{\Tor}{\operatorname{Tor}}

\renewcommand{\H}{\operatorname{H}}
\newcommand{\rH}{\widetilde{\operatorname{H}}}
\newcommand{\rC}{\widetilde{C}}

\renewcommand{\P}{\operatorname{P}}
\newcommand{\tensor}{\otimes}

\newcommand{\Ker}{\operatorname{Ker}}
\renewcommand{\Im}{\operatorname{Im}}

\newcommand{\shift}{s}
\newcommand{\dbl}{\mathopen{[\![}}
\newcommand{\dbr}{\mathopen{]\!]}}

\newcommand{\zn}{\numbers^t}
\newcommand{\cL}{cL_I}
\newcommand{\sat}[1]{K(#1)}
\newcommand{\satt}[1]{\hat{K}(#1)}
\newcommand{\sattt}[1]{\bar{K}(#1)}
\newcommand{\saturation}[1]{\bar{#1}}

\newcommand{\free}[1]{\Lambda #1}
\newcommand{\freep}[1]{\Lambda(#1)}

\newtheorem{theorem}{Theorem}
\newtheorem{thm}{Theorem}
\newtheorem{lemma}{Lemma}
\newtheorem{corollary}{Corollary}
\newtheorem{cor}{Corollary}

\theoremstyle{definition}
\newtheorem{remark}{Remark}

\begin{document}

\author{Alexander Berglund}
\title{Poincar\'e series of monomial rings}

\address{Department of Mathematics, Stockholm University, SE-106 91
Stockholm, Sweden}

\email{alexb@math.su.se}


\begin{abstract}
Let $k$ be a field, let $I$ be an ideal generated by monomials in the
polynomial ring $k[x_1,\ldots,x_t]$ and let $R=k[x_1,\ldots,x_t]/I$ be
the associated monomial ring. The $k$-vector spaces $\Tor_i^R(k,k)$
are $\zn$-graded. We derive a formula for the multigraded Poincar\'e
series of $R$,
$$\P_k^R(\boldx,z) = \sum_{i\geq 0, \alpha\in\zn}
\dim_k\Tor_{i,\alpha}^R(k,k)x^\alpha z^i,$$ in terms of the homology
of certain simplicial complexes associated to subsets of the minimal
set of generators for $I$. The homology groups occuring in the formula
can be interpreted as the homology groups of lower intervals in the
\emph{lattice of saturated subsets} of the generators for $I$.
\end{abstract}

\maketitle

\section{Introduction}
Let $I$ be an ideal generated by monomials in the polynomial ring $Q=
k[\boldx]=k[x_1,\ldots,x_t]$ over some field $k$ and let $R=Q/I$ be
the associated monomial ring. The $\zn$-grading of $Q$ assigning the
degree $\alpha = (\alpha_1,\ldots,\alpha_t)$ to the monomial $x^\alpha
= x_1^{\alpha_1}\cdot\ldots\cdot x_t^{\alpha_t}$ is inherited by $R$
and by $k\cong Q/(\boldx)$. Thus the $k$-vector space $\Tor_i^R(k,k)$
can be equipped with an $\zn$-grading for each $i$. The formal power
series
$$\P_k^R(\boldx,z) = \sum_{i\geq 0,\alpha\in\zn}
\dim_k\Tor_{i,\alpha}^R(k,k)x^\alpha z^i
\in\integers\dbl x_1,\ldots,x_t,z \dbr$$ is called the \emph{multigraded
Poincar\'e series} of $R$. It is proved by Backelin in \cite{backelin}
that this series is rational of the form
$$\P_k^R(\boldx,z) = \frac{\prod_{i=1}^t(1+x_iz)}{b_R(\boldx,z)},$$
for a polynomial $b_R(\boldx,z)\in\integers[x_1,\ldots,x_t,z]$.

In this paper we derive a formula for $b_R(\boldx,z)$ in terms of the
homology of certain simplicial complexes $\Delta_S'$ associated to
subsets $S$ of the minimal set of generators, $M_I$, for $I$.

Sets of monomials are considered as undirected graphs by letting edges
go between monomials having non-trivial common factors. Denote by
$\components{S}$ the number of connected components of a monomial set
$S$. If $S$ is a finite set of monomials let $m_S$ denote the least
common multiple of all elements of $S$.

Given a pair $S\subseteq M$ of monomial sets, $S$ is called
\emph{saturated} in $M$ if for all $m\in M$ and all connected subsets
$T$ of $S$, $m$ divides $m_T$ only if $m\in S$. Denote by $\sat{M}$
the set of non-empty saturated subsets of $M$.

If $M$ is a monomial set and $M=M_1\cup\ldots\cup M_r$ is its
decomposition into connected components, then $\Delta_M'$ is the
simplicial complex with vertex set $M$ and simplices
$$\set{S\subseteq M}{m_S\ne m_M\,\,\mbox{or $M_i\cap S$ disconnected
for some $i$}}.$$

If $H = \bigoplus_{i\in\integers}H_i$ is a graded vector space with
$H_i$ finite dimensional for each $i$, then denote $H(z) = \sum_{i\in
\integers} \dim H_i z^i \in \integers \dbl z^{-1},z \dbr$.

After these explanations we can state the main result.
\begin{theorem} \label{hauptsatz}
Let $k$ be any field. Let $I$ be an ideal in $Q=k[x_1,\ldots,x_t]$
generated by monomials of degree at least $2$, and let $M$ be its
minimal set of generators. The denominator of the Poincar\'e series of
$R=Q/I$ is given by
\begin{equation} \label{hauptformel}
b_R(\boldx,z) = 1+\sum_{S\in\sat{M}}
m_S(-z)^{\components{S}+2}\rH(\Delta_S';k)(z),
\end{equation}
\end{theorem}

In section \ref{applications} we will see that $\rH(\Delta_S';k)$ can
be interpreted as the reduced homology of the open interval
$(\emptyset,S)$ in the set $\sat{M}$ partially ordered by inclusion.

Theorem \ref{hauptsatz} has the following corollary, which gives a
solution to a problem posed by Avramov in \cite{avramov}. It improves
the bound of the $z$-degree of $b_R(\boldx,z)$ given in
\cite{backelin} considerably.

\begin{corollary} \label{degree bound}
With notations as in Theorem \ref{hauptsatz}
$$\deg b_R(z) \leq n+g,$$ where $b_R(z) = b_R(1,\ldots,1,z)$, $n=|M|$
is the number of minimal generators of $I$ and $g$ is the largest size
of a discrete subset of $M$. In particular
$$\deg b_R(z) \leq 2n,$$ with equality if and only if $R$ is a
complete intersection.
\end{corollary}

For a monomial set $M$, let $L_M$ denote the set $\set{m_S}{S\subseteq
M}$ partially ordered by divisibility. If $I$ is a monomial ideal, the
set $L_I = L_{M_I}$ is called the \emph{lcm-lattice} of $I$,
cf. \cite{gasharov}.

Let $I$ and $I'$ be monomial ideals in the polynomial rings
$Q=k[\boldx]$ and $Q'=k[\boldx']$, respectively, where $\boldx$ and
$\boldx'$ are finite sets of variables. One can introduce an
equivalence relation on the class of monomial ideals by declaring that
$I$ and $I'$ are equivalent if $L_I$ and $L_{I'}$ are isomorphic as
partially ordered graphs, that is, if there is a bijection $f\colon
L_I\rightarrow L_{I'}$ which respects both structures. Note that we do
not require that the ideals live in the same polynomial ring. This
equivalence relation is studied by Gasharov, Peeva and Welker in
\cite{gasharov}, where they prove that the ring $R=Q/I$ is Golod if
and only if $R'=Q'/I'$ is.

Avramov proves in \cite{avramov} that, with $R$ and $R'$ as above,
there is an isomorphism of graded Lie algebras
$$\pi^{\geq 2}(R) \cong \pi^{\geq 2}(R').$$ Here $\pi(R)$ is the
\emph{homotopy Lie algebra} of $R$, cf. \cite{avramov-inf} chapter
$10$. A corollary of this result is the equality
$$b_R(z) = b_{R'}(z),$$ where $b_R(z) = b_R(1,\ldots,1,z)$. Thus after
fixing the coefficient field, $b_R(z)$ depends only on the equivalence
class of $I$. This generalizes the fact that whether
$R$ is Golod or not depends only on the equivalence class of $I$.

As a consequence of our formula we obtain a slight strengthening of
Avramov's corollary, namely
\begin{corollary} \label{equivalence}
Let $I$ and $I'$ be ideals generated by monomials of degree at least
$2$ in the rings $k[\boldx]$ and $k[\boldx']$ respectively, where
$\boldx$ and $\boldx'$ are finite sets of variables and $k$ is a
field. Let $R=k[\boldx]/I$, $R' = k[\boldx']/I'$. If $f:L_I\rightarrow
L_{I'}$ is an isomorphism of partially ordered graphs, then
$$b_{R'}(\boldx',z) = f(b_R(\boldx,z)),$$ where $f(b_R(\boldx,z))$
denotes the result of applying $f$ to the coefficients of
$b_R(\boldx,z)$, regarding it as a polynomial in $z$.
\end{corollary}

Theorem \ref{hauptsatz} is first proved in the case when $I$ is
generated by square-free monomials and then a construction by
Fr\"oberg \cite{froberg} is used to reduce the general case to this
case. The starting point of the proof in the square-free case is the
observation that the Poincar\'e series $\P_k^R(\boldx,z)$ is
determined by a finite set of the \emph{multigraded deviations}
$\epsilon_{i,\alpha}(R)$ of the ring $R$. Then we use the fact that
these multigraded deviations can be computed from a \emph{minimal
model} for $R$ over the polynomial ring.

Section \ref{combinatorics} contains definitions and conventions
concerning relevant combinatorial notions. The facts about minimal
models needed are given in section \ref{models}. The proof of Theorem
\ref{hauptsatz}, along with some auxiliary results, is presented in
section \ref{proof section}. Section \ref{applications} contains the
corollaries of the main result and their proofs. In the concluding
remark of section \ref{applications}, we note how Theorem
\ref{hauptsatz} gives a combinatorial criterion for a monomial ring to
be Golod.

\section{Combinatorics} \label{combinatorics}

\subsection{Simplicial complexes.}
A simplicial complex on a set $V$ is a set $\Delta$ of subsets of $V$
such that $F\subseteq G\in \Delta$ implies $F\in \Delta$. $V$ is the
\emph{vertex set} of $\Delta$. The \emph{$i$-faces} or
\emph{$i$-simplices} of $\Delta$ are the elements in $\Delta$ of
cardinality $i+1$. We do not require that $\{v\}\in \Delta$ for all
$v\in V$, but if a simplicial complex $\Delta$ is given without
reference to a vertex set $V$, then it is assumed that $V=\cup
\Delta$.

If $\Delta$ is a simplicial complex then $\rC(\Delta)$ will denote the
augmented chain complex associated to $\Delta$. Thus $\rC_i(\Delta)$
is the free abelian group on the $i$-faces of $\Delta$, $\emptyset$
being considered as the unique $(-1)$-face, and $\rC(\Delta)$ is
equipped with the standard differential of degree $-1$. Therefore
$$\H_i(\rC(\Delta)) = \rH_i(\Delta).$$ As usual, if $G$ is an abelian
group, then $\rC(\Delta;G) = \rC(\Delta)\tensor G$ and
$\rH_i(\Delta;G) = \H_i(\rC(\Delta;G))$.

The Alexander dual of a simplicial complex $\Delta$ with vertices $V$
is the complex
$$\Delta^\vee = \set{F\subseteq V}{V\setminus F\not\in\Delta}.$$ The
join of two complexes $\Delta_1$, $\Delta_2$ with disjoint vertex sets
$V_1,V_2$ is the complex with vertex set $V_1\cup V_2$ and faces
$$\Delta_1*\Delta_2 = \set{F_1\cup F_2}{F_1\in\Delta_1,
F_2\in\Delta_2}.$$ With $\Delta_1$ and $\Delta_2$ as above, we define
what could be called the \emph{dual join} of them:
$$\Delta_1\djoin\Delta_2 = (\Delta_1^\vee * \Delta_2^\vee)^\vee.$$ Thus
$\Delta_1\djoin\Delta_2$ is the simplicial complex with vertex set $V_1\cup
V_2$ and simplices
$$\set{F\subseteq V_1\cup V_2}{F\cap V_1\in\Delta_1\,\,\mbox{or}\,\,
F\cap V_2 \in\Delta_2}.$$

We will now briefly review the effects of these operations on the
homology groups when the coefficients come from a field $k$.

If $|V| = n$, then (\cite{bruns} Lemma $5.5.3$)
$$\rH_i(\Delta;k) \cong \rH_{n-i-3}(\Delta^\vee;k).$$ If $C$ is a
chain complex, then $\shift C$ denotes the shift of $C$; $(\shift C)_i
= C_{i-1}$. Because of the convention that a set with $d$ elements has
dimension $d-1$ considered as a simplex there is a shift in the
following formula:
$$\rH(\Delta_1*\Delta_2;k) \cong
\shift(\rH(\Delta_1;k)\tensor_k\rH(\Delta_2;k)).$$

If $H = \bigoplus_{i\in\integers} H_i$ is a graded vector space such
that each $H_i$ is of finite dimension, then let $H(z) =
\sum_{i\in\integers} \dim H_i z^i$ be the \emph{generating function}
of $H$. The above isomorphisms of graded vector spaces can be
interpreted in terms of generating functions. If $\Delta$ has $n$
vertices, then
\begin{equation*}
z^n\rH(\Delta^\vee;k)(z^{-1}) = z^3\rH(\Delta;k)(z),
\end{equation*}
and if $\Delta = \Delta_1*\Delta_2$, then
\begin{equation*}
\rH(\Delta;k)(z) = z\rH(\Delta_1;k)(z)\cdot\rH(\Delta_2;k)(z).
\end{equation*}
From these two identities and an induction one can work out the
following formula. If $\Delta = \Delta_1\djoin\ldots\djoin\Delta_r$,
then
\begin{equation} \label{sum}
\rH(\Delta;k)(z) =
z^{2r-2}\rH(\Delta_1;k)(z)\cdot\ldots\cdot\rH(\Delta_r;k)(z).
\end{equation}

\subsection{Sets of monomials.} \label{monomials}
Let $x_1,\ldots,x_t$ be variables. If $\alpha\in\numbers^t$, then we
write $x^\alpha$ for the monomial $ x_1^{\alpha_1}\cdot\ldots\cdot
x_t^{\alpha_t}$. The \emph{multidegree} of $x^\alpha$ is $\mdeg
x^\alpha = \alpha$. If $\alpha\in\{0,1\}^t$, then both $\alpha$ and
$x^\alpha$ are called \emph{square-free}.

To a set $M$ of monomials we associate an undirected graph, with
vertices $M$, whose edges go between monomials having a non-trivial
common factor. This is the graph structure referred to when properties
such as connectedness et.c., are attributed to monomial sets. Thus,
for instance, a monomial set is called \emph{discrete} if the
monomials therein are pairwise without common factors. By $D(M)$ we
denote the set of non-empty discrete subsets of $M$. A connected
component of $M$ is a maximal connected subset. Any monomial set $M$
has a decomposition into connected components $M= M_1 \cup \ldots \cup
M_r$, and we let $\components{M} = r$ denote the number of such.

If $I$ is an ideal in a polynomial ring generated by monomials there
is a uniquely determined minimal set of monomials generating $I$.
This minimal generating set, denoted $M_I$, is characterized by being
an \emph{antichain}, that is, for all $m,n\in M_I$, $m|n$ implies
$m=n$.

If $S$ is a finite set of monomials, then $m_S$ denotes the least
common multiple of all elements of $S$. By convention $m_\emptyset =
1$. The set $L_M = \set{m_S}{S\subseteq M}$ partially ordered by
divisibility is a lattice with $\lcm$ as join, called the
\emph{lcm-lattice} of the set $M$. If $I$ is a monomial ideal, then
$L_I:=L_{M_I}$ is called the lcm-lattice of $I$.

Remark: The \emph{gcd-graph} of $I$, studied in \cite{avramov}, is the
complement of the graph $L_I$.

If $M,N$ are two sets of monomials then $M_N$ denotes the set of those
monomials in $M$ which divide some monomial in $N$. Write $M_m =
M_{\{m\}}$, and $M_\alpha = M_{x^\alpha}$.

Let $S$ be a subset of a monomial set $M$. The \emph{saturation} of
$S$ in $M$ is the set $\saturation{S} = M_N$, where $N =
\{m_{S_1},\ldots,m_{S_r}\}$ if $S_1,\ldots,S_r$ are the connected
components of $S$. Clearly $S\subseteq\saturation{S}$, and $S$ is
called \emph{saturated in $M$} if equality holds. Equivalently, $S$ is
saturated in $M$ if for all $m\in M$, $m\mid m_T$ implies $m\in S$ if
$T$ is a connected subset of $S$. Clearly $S$ is saturated in $M$ if
and only if all the connected components of $S$ are. The set of
saturated subsets of $M$ is denoted $\satt{M}$, and the set of
non-empty such subsets is denoted $\sat{M}$.

Two monomial sets $M$, $N$ are called \emph{equivalent} if there is an
isomorphism of partially ordered sets $f\colon L_M\rightarrow L_N$
which is also an isomorphism of graphs. Such a map $f$ will be called
an \emph{equivalence}.

To a monomial set $M$ we associate a simplicial complex $\Delta_M$,
with vertex set $M$ and faces
$$\set{S\subseteq M}{m_S\ne m_M\,\,\mbox{or $S$ disconnected}}.$$

We state again the definition of the simplicial complex $\Delta_M'$
given in the introduction. Let $M=M_1\cup\ldots\cup M_r$ be the
decomposition of $M$ into its connected components. The complex
$\Delta_M'$ has vertex set $M$ and faces
$$\set{S\subseteq M}{m_S\ne m_M\,\,\mbox{or $M_i\cap S$ disconnected
for some $i$}}.$$ Note that
$$\Delta_M' = \Delta_{M_1}\djoin\ldots\djoin\Delta_{M_r}.$$

\begin{lemma} \label{equivalent-antichains}
Let $M$ and $N$ be antichains of monomials. An equivalence $f\colon
L_M\rightarrow L_N$ induces a bijection $\sat{M}\rightarrow \sat{N}$,
where $S\in\sat{M}$ is mapped to $f(S)\in\sat{N}$. Furthermore, for
every $S\subseteq M$, $S$ and $f(S)$ are isomorphic as graphs and the
complexes $\Delta_S'$, $\Delta_{f(S)}'$ are isomorphic.
\end{lemma}

\begin{proof}
If $M$ is an antichain, then $L_M$ is atomic with atoms $M$. An
isomorphism of lattices maps atoms to atoms, so $f$ restricts to a
graph isomorphism $M\rightarrow N$ and hence gives rise to a bijection
of subgraphs of $M$ to isomorphic subgraphs of $N$, which clearly maps
saturated sets to saturated sets. Since the definition of $\Delta_S'$
is phrased in terms of the graph structure of $S$ and on the lattice
structure of $L_S\subseteq L_M$, it is clear that
$\Delta_S'\cong\Delta_{f(S)}'$ for any subset $S$ of $M$.
\end{proof}

If $k$ is a field and $M$ is a monomial set in the variables
$x_1,\ldots,x_t$ then $R = k[x_1,\ldots,x_t]/(M)$ is the monomial ring
associated to $M$. Monomial ideals are homogeneous with respect to the
multigrading of $k[x_1,\ldots,x_t]$, so monomial rings inherit this
grading.

\section{$\zn$-graded models and deviations} \label{models}
In this section we will collect and adapt to the $\zn$-graded
situation some well known results on models of commutative
rings. There is no claim of originality. Our main reference is
\cite{avramov-inf} chapter $7.2$. Inspiration comes also from the
sources \cite{halperin} and \cite{gulliksen-levin}, where analogous
but not directly applicable results can be found.

The notation $|x|$ refers to the homological degree of an element
$x$. We use $\mdeg x$ to denote the multidegree of $x$.

\subsection{Deviations.}
Let $R = Q/I$, where $I$ is a monomial ideal in $Q$. Recall that the
$i$:th \emph{deviation}, $\epsilon_i = \epsilon_i(R)$, of the ring $R$
can be defined as the number of variables adjoined in degree $i$ in an
acyclic closure, $R\langle X \rangle$, of $k$ over $R$,
cf. \cite{avramov-inf} Theorem $7.1.3$. See \cite{avramov-inf} section
$6.3$ for the construction of acyclic closures. The ring $R$ is
$\zn$-graded and $k$ is an $\zn$-graded $R$-module, and one can show
that there is a unique $\zn$-grading on the acyclic closure $R\langle
X \rangle$ which is respected by the differential. One may therefore
introduce $\zn$-graded deviations
$$\epsilon_{i,\alpha}=|\set{x\in X}{|x|=i,\,\,\mdeg{x} = \alpha}|.$$
It is clear that $\epsilon_{i,\alpha}=0$ if $|\alpha|<i$.

By the general theory, $R\langle X \rangle$ is a minimal resolution of
$k$, and hence $R\langle X \rangle\tensor_R k = \\ \H(R\langle X
\rangle\tensor_R k)$ is isomorphic as a multigraded vector space to
$\Tor^R(k,k)$. This yields a product representation of the multigraded
Poincar\'e series
\begin{eqnarray} \label{product}
\P_k^R(\boldx,z) & = & \prod_{i\geq 1,\alpha\in\zn}\frac{(1+x^\alpha
z^{2i-1})^{\epsilon_{2i-1,\alpha}}}{(1-x^\alpha
z^{2i})^{\epsilon_{2i,\alpha}}}.
\end{eqnarray}

It is a fundamental result that the deviations $\epsilon_{i}$ can be
computed from a \emph{minimal model} of $R$ over $Q$,
cf. \cite{avramov-inf} Theorem $7.2.6$. Only trivial modifications are
required in order to show that the $\zn$-graded deviations
$\epsilon_{i,\alpha}$ can be computed from an $\zn$-graded minimal
model of $R$ over $Q$. For this reason, we will need a few facts about
minimal models.

\subsection{Free dg-algebras.}
Let $V=\bigoplus_{i\geq 0} V_i$ be a graded vector space over $k$ such
that $\dim_k V_i < \infty$ for each $i$. We denote by $\free{V}$ the
free graded commutative algebra on $V$, that is,
$$\free{V} = \mbox{exterior algebra}(V_{\operatorname{odd}}) \tensor_k
\mbox{symmetric algebra}(V_{\operatorname{even}}).$$

Denote by $(V)$ the ideal generated by $V$ in $\free{V}$. A
homomorphism $f\colon\free{V}\rightarrow\free{W}$ of graded algebras
with $f(V)\subseteq (W)$ induces a linear map $Lf\colon V\rightarrow
W$, called the \emph{linear part} of $f$, which is defined by the
requirement $f(v)-Lf(v)\in (W)^2$ for all $v\in V$.

If $x_1,\ldots,x_t$ is a basis for $V_0$, then $\free{V} =
Q\tensor_k\freep{V_+}$, where $Q=k[x_1,\ldots,x_t]$ and $V_+$ is the
sum of all $V_i$ for positive $i$. Therefore $\free{V}$ may be
regarded as a $Q$-module and each $(\free{V})_n$ is a finitely
generated free $Q$-module. Let $\mathfrak{m}\subseteq Q$ be the
maximal ideal generated by $V_0$ in $Q$. Note that $(V) =
(V_+)+\mathfrak{m}$ as vector spaces. The following basic lemma is a
weak counterpart of Lemma $14.7$ in \cite{halperin} and of Lemma
$1.8.7$ in \cite{gulliksen-levin}.

\begin{lemma} \label{algebra isomorphism}
Let $f\colon\free{U}\rightarrow\free{V}$ be a homomorphism of graded
algebras such that $f_0\colon\freep{U_0}\rightarrow\freep{V_0}$ is an
isomorphism and the linear part, $Lf\colon U\rightarrow V$, is an
isomorphism of graded vector spaces. Then $f$ is an isomorphism.
\end{lemma}

\begin{proof}
Identify $Q=\freep{U_0}=\freep{V_0}$ via $f_0$. Since $Lf$ is an
isomorphism, $\free{U}$ and $\free{V}$ are isomorphic. Thus to show
that $f$ is an isomorphism it is enough to show that $f_n\colon
(\free{U})_n\rightarrow (\free{V})_n$ is surjective in each degree
$n$, because $f_n$ is a map between finitely generated isomorphic
$Q$-modules. We do this by induction. The map $f_0$ is surjective by
assumption. Let $n\geq 1$ and assume that $f_i$ is surjective for
every $i<n$. Then since $Lf$ is surjective we have
$$(\free{V})_n\subseteq f((\free{U})_n) + ((V_+)^2)_n +
\mathfrak{m}(\free{V})_n.$$ $((V_+)^2)_n$ is generated by products
$vw$, where $|w|,|v|<n$, so by induction $((V_+)^2)_n\subseteq
f((\free{U})_n)$. Hence
$$(\free{V})_n\subseteq f((\free{U})_n) + \mathfrak{m}(\free{V})_n.$$
$(\free{V})_n$ and $f((\free{U})_n)$ are graded $Q$-modules, so it
follows from the graded version of Nakayama's lemma that $(\free{V})_n
= f((\free{U})_n)$.
\end{proof}

By a \emph{free dg-algebra}, we will mean a dg-algebra of the form
$(\free{V},d)$, for some graded vector space $V$, where $d$ is a
differential of degree $-1$ satisfying $dV\subseteq (V)$. The linear
part $Ld$ of the differential $d$ on $\free{V}$ is a differential on
$V$, and will be denoted $d_0$. A free dg-algebra $(\free{V},d)$ is
called \emph{minimal} if $dV\subseteq (V)^2$. Thus $(\free{V},d)$
is minimal if and only if $d_0=0$.

If $(\free{V},d)$ is a free dg-algebra which is $\zn$-graded, that is,
there is a decomposition
$$(\free{V})_i = \bigoplus_{\alpha\in\zn} (\free{V})_{i,\alpha}$$ such
that $d(\free{V})_{i,\alpha}\subseteq (\free{V})_{i-1,\alpha}$, then
denote
$$\H_{i,\alpha}(\free{V},d) = \H_i((\free{V})_\alpha,d).$$ The
$\zn$-grading is called non-trivial if $\mdeg v\ne 0$ for all $v\in V$.

The following is a counterpart of Lemma $3.2.1$ in
\cite{gulliksen-levin}, but taking the $\zn$-grading into account. It
tells us how to `minimize' a given dg-algebra.

\begin{lemma} \label{minimal}
Let $(\free{V},d)$ be a non-trivially $\zn$-graded dg-algebra with
$dV_1\subseteq \mathfrak{m}^2$. Then there exists an $\zn$-graded
minimal dg-algebra $(\free{H},d')$ with $H\cong\H(V,d_0)$, together
with a surjective morphism of dg-algebras
$$(\free{V},d)\rightarrow (\free{H},d')$$ which is a quasi-isomorphism
if $k$ has characteristic $0$. For arbitrary $k$ we have
$\H_0(\free{V},d)\cong\H_0(\free{H},d')$ and
$$\H_{i,\alpha}(\free{V},d) \cong \H_{i,\alpha}(\free{H},d')$$ for all
square-free $\alpha\in\zn$ and all $i$.
\end{lemma}

\begin{proof}
Let $W$ be a graded subspace of $V$ such that $V=\Ker d_0 \oplus W$
and similarly split $\Ker d_0$ as $H\oplus \Im d_0$ (hence
$H\cong\H(V,d_0)$). Note that since $dV_1\subseteq\mathfrak{m}^2$,
$W_0=W_1=0$. $d_0$ induces an isomorphism $W\rightarrow \Im d_0$, so
we may write
$$V = H\oplus W\oplus d_0(W).$$ Consider the graded subspace $U =
H\oplus W \oplus dW$ of $\free{V}$. The induced homomorphism of graded
algebras
$$f\colon \free{U}\rightarrow \free{V}$$ is an isomorphism by Lemma
\ref{algebra isomorphism}, because $f_0$ is the identity on
$\free{H_0}$ and the linear part of $f$ is the map $1_H\oplus
1_W\oplus g$, where $g\colon dW\rightarrow d_0(W)$ is the isomorphism
taking an element to its linear part (isomorphism precisely because
$\Ker d_0\cap W=0$). Thus we may identify $\free{U}$ and $\free{V}$
via $f$. In particular $f^{-1}df$ is a differential on $\free{U}$,
which we also will denote by $d$, and $(\free{U},d)$ is a dg-algebra
in which $\freep{W\oplus dW}$ is a dg-subalgebra.

The projection $U\rightarrow H$ induces an epimorphism of graded
algebras
$$\phi\colon \free{U}\rightarrow\free{H}$$ with kernel $(W\oplus
dW)\free{U}$, the ideal generated by $W\oplus dW$ in
$\free{U}$. Define a differential $d'$ on $\free{H}$ by
$$d'(h) = \phi d\iota(h),$$ where $\iota$ is induced by the inclusion
$H\subseteq U$. With this definition $\phi$ becomes a morphism of
dg-algebras and it is evident that $(\free{H},d')$ is
minimal. Furthermore $\H(H,d_0')=H\cong \H(V,d_0)$ by definition.

Consider the increasing filtration
$$F_p = (\free{H})_{\leq p}\cdot\freep{W\oplus dW}.$$ Obviously $\cup
F_p = \free{U}$, and $dF_p \subseteq F_p$ since $d$ preserves
$\freep{W\oplus dW}$. The associated first quadrant spectral sequence
is convergent, with
$$E^2_{p,q} = \H_p(\free{H},d')\tensor_k\H_q(\freep{W\oplus dW},d)
\Longrightarrow \H_{p+q}(\free{U},d).$$ Since $W_0 = W_1 = 0$, we have
$\H_0(\freep{W\oplus dW},d) = k$, and therefore $\H_0(\free{H},d') =
E^2_{0,0} = E^3_{0,0} = \ldots = E^\infty_{0,0} = \H_0(\free{U},d) =
\H_0(\free{V},d)$. If the field $k$ has characteristic zero, then
$(\freep{W\oplus dW},d)$ is acyclic, so in this case the spectral
sequence collapses, showing that
$\H(\free{H},d')\cong\H(\free{V},d)$. However, $\freep{W\oplus dW}$
need not be acyclic in positive characteristic $p$ --- if $x\in W$ is
of even degree, then $x^{np}$ and $x^{np-1}dx$ represent non-trivial
homology classes for all $n\geq 1$. Recall however that we are working
with $\zn$-graded objects and maps. Since we have a non-trivial
$\zn$-grading, $\freep{W\oplus dW}_\alpha$ is acyclic for square-free
$\alpha$, simply because no elements of the form $x^na$, for $x\in
(W\oplus dW),a\in\freep{W\oplus dW},n>1$, are there. In particular the
dissidents $x^{np}$ and $x^{np-1}dx$ live in non-square-free
degrees. Hence the spectral sequence collapses in square-free degrees,
regardless of characteristic, and so
$$\H_{i,\alpha}(\free{H},d')\cong\H_{i,\alpha}(\free{V},d),$$ for all
square-free $\alpha$ and all $i$.
\end{proof}

\subsection{Models.}
Let $Q=k[x_1,\ldots,x_t]$ and let $R = Q/I$ for some ideal $I\subseteq
\mathfrak{m}^2$. A \emph{model} for the ring $R$ over $Q$ is a free
dg-algebra $(\free{V},d)$ with $(\free{V})_0 = Q$, such that
$\H_0(\free{V},d) = R$ and $\H_i(\free{V},d) = 0$ for all $i>0$. In
particular a model for $R$ over $Q$ is a free resolution
$$\cdots \rightarrow (\free{V})_n \rightarrow
(\free{V})_{n-1}\rightarrow \cdots \rightarrow (\free{V})_1\rightarrow
Q \rightarrow R\rightarrow 0$$ of $R$ as a $Q$-module. A model
$(\free{V},d)$ is called \emph{minimal} if it is a minimal
dg-algebra. A minimal model for $R$ over $Q$ always exists, and is
unique up to (non-canonical) isomorphism, cf. \cite{avramov-inf}
Proposition $7.2.4$.

If $R$ is $\zn$-graded, one can ask for the minimal model of $R$ to be
$\zn$-graded.

\begin{lemma} \label{build model}
Let $(\free{V},d)$ be a minimal $\zn$-graded dg-algebra with
$\H_0(\free{V},d) = R$, and assume that
$$\H_{i,\alpha}(\free{V},d) = 0$$ for all $i>0$ and all square-free
$\alpha\in \zn$. Then $(\free{V},d)$ can be completed to a minimal
model $(\free{W},d)$ of $R$ such that $W_\alpha = V_\alpha$ for all
square-free $\alpha$.
\end{lemma}

\begin{proof}
A minimal model can be constructed inductively, by successively
adjoining basis elements to $V$ in order to kill homology,
cf. \cite{avramov-inf} propositions $2.1.10$ and $7.2.4$ for
details. Since $\free{V}$ is $\zn$-graded, we can do this one
multidegree at a time. Adding a basis element of multidegree $\alpha$
will not affect the part of the algebra below $\alpha$. Since
$\H_i((\free{V})_\alpha) = 0$ for all $i>0$ when $\alpha$ is
square-free, we do not need to add variables of square-free
multidegrees in order to kill homology. Thus, applying this technique,
we get a minimal model $\free{W}$ of $R$, where $W$ is a vector space
obtained from $V$ by adjoining basis elements of non-square-free
degrees. In particular $W_\alpha = V_\alpha$ for all square-free
$\alpha$.
\end{proof}

Taking $\zn$-degrees into account, it is not difficult to modify the
proof of Theorem $7.2.6$ in \cite{avramov-inf} to obtain the following
result:
\begin{lemma} \label{deviations}
Let $(\free{W},d)$ be an $\zn$-graded minimal model for $R$ over
$Q$. Then the $\zn$-graded deviations $\epsilon_{i,\alpha}$ of $R$ are
given by
$$\epsilon_{i,\alpha} = \dim_k W_{i-1,\alpha},$$ for $i\geq
1, \alpha\in\zn$. \hfill $\square$
\end{lemma}

\section{Poincar\'e series} \label{proof section}
This section is devoted to the deduction of Theorem \ref{hauptsatz}.
\begin{thm}
Let $k$ be any field. Let $I$ be an ideal in $Q=k[x_1,\ldots,x_t]$
generated by monomials of degree at least $2$, and let $M$ be its
minimal set of generators. The denominator of the Poincar\'e series of
$R=Q/I$ is given by
$$b_R(\boldx,z) = 1+\sum_{S\in\sat{M}}
m_S(-z)^{\components{S}+2}\rH(\Delta_S';k)(z).$$
\end{thm}
Some intermediate results will be needed before we can give the
proof. Retain the notations of the theorem throughout this section. We
will frequently suppress the variables and write $b_R = b_R(\boldx,z)$
and $\P_k^R = \P_k^R(\boldx,z)$.

\subsection{An observation.}
Assume to begin with that the ideal $I$ is minimally generated by
square-free monomials $M=\{m_1,\ldots,m_n\}$ of degree at least
$2$. If we are given a subset $S=\{m_{i_1},\ldots,m_{i_r}\}$ of $M$,
where $i_1<\ldots<i_r$, then set $\sgn{m_{i_j},S} = (-1)^{j-1}$. By
Backelin \cite{backelin}, the Poincar\'e series of $R$ is rational of
the form
$$\P_k^R(\boldx,z) = \frac{\prod_{i=1}^t(1+x_iz)}{b_R(\boldx,z)},$$
where $b_R(\boldx,z)$ is a polynomial with integer coefficients and
$x_i$-degree at most $1$ for each $i$. We start with the following
observation made while scrutinizing Backelin's proof.

\begin{lemma} \label{observation}
If $I$ is generated by square-free monomials, then the polynomial
$b_R$ is square-free with respect to the $x_i$-variables. Moreover
$b_R$ depends only on the deviations $\epsilon_{i,\alpha}$ for
square-free $\alpha$. In fact, there is a congruence modulo
$(x_1^2,\ldots,x_t^2)$:
$$b_R \equiv \prod_{\alpha\in \{0,1\}^t} (1-x^\alpha p_\alpha(z)),$$
where $p_\alpha(z)$ is the polynomial $p_\alpha(z) =
\sum_{i=1}^{|\alpha|} \epsilon_{i,\alpha}z^i$.
\end{lemma}

\begin{proof}
Note that $\epsilon_{1,e_i} = 1$ and $\epsilon_{1,\alpha} = 0$ for
$\alpha\ne e_i$ ($i=1,\ldots,t$). Hence using the product
representation \eqref{product} and reducing modulo
$(x_1^2,\ldots,x_t^2)$ we get (note that $(1+mp(z))^n\equiv 1+nmp(z)$
for any integer $n$ and any square-free monomial $m$):
\begin{eqnarray} \label{formula}
b_R & = & \frac{\prod_{i\geq 1,\alpha}(1-x^\alpha
    z^{2i})^{\epsilon_{2i,\alpha}}}{\prod_{i\geq 2,\alpha}(1+x^\alpha
    z^{2i-1})^{\epsilon_{2i-1,\alpha}}}\nonumber \\ & \equiv &
    \prod(1-x^\alpha
    (\epsilon_{2i-1,\alpha}z^{2i-1}+\epsilon_{2i,\alpha}z^{2i}))
    \nonumber \\ & \equiv & \prod(1-x^\alpha p_{\alpha}(z)),
\end{eqnarray}
product taken over all square-free $\alpha$, where
$p_\alpha(z)\in\integers[z]$ is the polynomial $p_\alpha(z) =
\sum_{i=1}^{|\alpha|}\epsilon_{i,\alpha}z^i$.
\end{proof}

This gives a formula for $b_R$ in terms of the square-free deviations
$\epsilon_{i,\alpha}$. Therefore we are interested in the square-free
part of an $\zn$-graded minimal model of $R$ over $Q$.

\subsection{Square-free deviations.}
In the square-free case, there is a nice interpretation of the
square-free deviations in terms of simplicial homology. Recall the
definition of $\Delta_M$ found in section \ref{monomials}. $M_\alpha$
denotes the set of monomials in $M$ which divide $x^\alpha$.

\begin{theorem} \label{squarefree deviations}
Assume that $I$ is minimally generated by a set $M$ of square-free
monomials of degree at least $2$. Let $\alpha$ be square-free and let
$i\geq 2$. If $x^\alpha\not\in L_I$, then $\epsilon_{i,\alpha} = 0$,
and if $x^\alpha\in L_I$ then
$$\epsilon_{i,\alpha} = \dim_k\rH_{i-3}(\Delta_{M_\alpha};k).$$
\end{theorem}

The proof of this theorem depends on the construction of the
square-free part of a minimal model for $R$, which we now will carry
out.

Let $C$ be the set of connected non-empty subsets of $M$ and let $V$
be the $\numbers\times\zn$-graded vector space with basis
$Y\cup\{x_1,\ldots,x_t;|x_i|=0,\mdeg{x_i} = e_i\}$, where
$$Y=\set{y_S}{S\in C, |y_S| = |S|, \mdeg{y_S} = \mdeg{m_S}}.$$ If $S$
is any subset of $M$ and $S=S_1\cup\ldots\cup S_r$ is its
decomposition into connected components, then define the symbol
$y_S\in\free{V}$ by
$$y_S = y_{S_1}\cdot\ldots\cdot y_{S_r}.$$ It follows at once that
$|y_S|=|S|$ and that $\mdeg{y_S} = \mdeg{m_S}$.

The differential $d$ on $\free{V}$ is defined on the basis $y_S$,
$S\in C$, by
\begin{equation} \label{diff}
dy_S = \sum_{s\in S}
\sgn{s,S}\frac{m_S}{m_{S\setminus\{s\}}}y_{S\setminus\{s\}},
\end{equation}
and is extended to all of $\free{V}$ by linearity and the Leibniz rule
(and of course $dx_i = 0$). Note that it may happen that
$y_{S\setminus\{s\}}$ becomes decomposable as a product in the sum
above. One verifies easily that the formula \eqref{diff} remains valid
for disconnected $S$. By definition $d$ is of degree $-1$ and respects
the $\zn$-grading.

Clearly, we have $\H_0(\free{V}) = R$. Let $\alpha\in\zn$ be
square-free. Then the complex $(\free{V})_\alpha$ is isomorphic to the
degree $\alpha$-part of the Taylor complex on the monomials $M$
(cf. \cite{froberg2}). It is well-known that the Taylor complex is a
resolution of $R$ over $Q$, so in particular $\H_i((\free{V})_\alpha)
= 0$ for all $i>0$. Since we assumed that the monomials $m_i$ are of
degree at least $2$, we have $dV_1 \subseteq
\mathfrak{m}^2$. Therefore, by Lemma \ref{minimal}, there is a minimal
$\zn$-graded dg-algebra $(\free{H},d')$ such that $H\cong
\H(V,d_0)$, $\H_0(\free{H},d') = \H_0(\free{V},d) = R$ and
$$\H_{i,\alpha}(\free{H},d') = \H_{i,\alpha}(\free{V},d) = 0,$$ for
all $i>0$ and all square-free $\alpha\in\zn$. Now, by Lemma \ref{build
model}, we can construct an $\zn$-graded minimal model $(\free{W},d)$
of $R$, such that $W_\alpha = H_\alpha$ for all square-free
$\alpha$. This is all we need to know about the minimal model
$(\free{W},d)$ in order to be able to prove Theorem \ref{squarefree
deviations}:

\begin{proof}[Proof of Theorem \ref{squarefree deviations}]
By Lemma \ref{deviations} we get that
\begin{equation} \label{deveq}
\epsilon_{i,\alpha} = \dim_k W_{i-1,\alpha} = \dim_k H_{i-1,\alpha} =
\dim_k \H_{i-1,\alpha}(V,d_0),
\end{equation}
for square-free $\alpha$. We will now proceed to give a combinatorial
description of the complex $V=(V,d_0)$. $V$ splits as a complex into
its $\zn$-graded components
$$V = \bigoplus_{\alpha\in\zn} V_\alpha.$$ $V_{e_i}$ is
one-dimensional and concentrated in degree $0$ for
$i=1,\ldots,t$. This accounts for the first deviations
$\epsilon_{1,e_i} = 1$. If $|\alpha|>1$, then $V_\alpha$ has basis
$y_S$ for $S$ in the set
$$C_\alpha=\set{S\subseteq M}{m_S =
x^\alpha,\,S\,\,\mbox{connected}}.$$ In particular $V_\alpha = 0$ if
$x^\alpha\not\in L_I$. The differential of $V_\alpha$ is given by
\begin{equation} \label{alpha-diff}
dy_S = \sum_{\substack{s\in S \\ S\setminus\{s\} \in C_\alpha}}
\sgn{s,S}y_{S\setminus\{s\}}.
\end{equation}

Let $\Sigma_\alpha$ be the simplicial complex whose faces are all
subsets of the set $M_\alpha = \set{m\in M}{m\mid x^\alpha}$, with
orientation induced from the orientation $\{m_1,\ldots,m_n\}$ of $M$.
Define a map from the chain complex $\rC(\Sigma_\alpha;k)$ to the
desuspended complex $s^{-1}V_\alpha$ by sending a face $S\subseteq
M_\alpha$ to $s^{-1}y_S$ if $S\in C_\alpha$ and to $0$ otherwise. In
view of \eqref{alpha-diff}, this defines a morphism of complexes,
which clearly is surjective. The kernel of this morphism is the chain
complex associated to $\Delta_{M_\alpha}$, so we get a short exact
sequence of complexes
$$0\rightarrow \rC(\Delta_{M_\alpha};k)\rightarrow
\rC(\Sigma_\alpha;k) \rightarrow s^{-1} V_\alpha \rightarrow 0$$ Since
$\Sigma_\alpha$ is acyclic, the long exact sequence in homology
derived from the above sequence shows that $\H_i(V_\alpha) \cong
\rH_{i-2}(\Delta_{M_\alpha};k)$. The theorem now follows from
\eqref{deveq}.
\end{proof}

In terms of the polynomials $p_\alpha(z)$ the theorem may be stated as
\begin{equation} \label{p_a}
p_\alpha(z) = z^3\rH(\Delta_{M_\alpha};k)(z),
\end{equation}
for $x^\alpha\in L_I$.

\subsection{Proof of Theorem \ref{hauptsatz}.}
\begin{proof}[Proof of Theorem \ref{hauptsatz}. Square-free case]
By Theorem \ref{squarefree deviations}, $p_\alpha(z)=0$ unless
$x^\alpha\in L_I$, in which case $p_\alpha(z) =
z^3\rH(\Delta_{M_\alpha};k)(z)$. But $\Delta_{M_\alpha}$ is
contractible if $M_\alpha$ is disconnected, so $p_\alpha(z) = 0$
unless $x^\alpha\in\cL$, where $\cL$ denotes the subset of $L_I$
consisting of elements $l\ne 1$ such that $M_l$ is connected.

Hence by Lemma \ref{observation}
$$b_R \equiv \prod_{x^\alpha\in\cL} (1-x^\alpha p_\alpha(z))\quad\mod
(x_1^2,\ldots,x_t^2).$$

If we carry out the multiplication in the above formula and use that
$b_R$ is square-free with respect to the $x_i$-variables (by Lemma
\ref{observation}) we get the equality
$$b_{R} = 1+\sum_{N\in D(\cL)} \prod_{x^\alpha\in
N}(-x^\alpha p_\alpha(z)) = 1+ \sum_{N\in D(\cL)}
m_N(-1)^{|N|}\prod_{x^\alpha\in N}p_\alpha(z)$$ (the identity
$\prod_{x^\alpha\in N} x^\alpha = m_N$ holds because $N$ is
discrete). Using \eqref{p_a} the formula takes the form
$$b_R = 1+\sum_{N\in D(\cL)}
m_N(-1)^{|N|}\prod_{x^\alpha\in N}
z^3\rH(\Delta_{M_\alpha};k)(z).$$ By \eqref{sum} this may
be written
$$b_R = 1+\sum_{N\in D(\cL)} m_N(-1)^{|N|} z^{|N|+2}
\rH(\Gamma;k)(z),$$ where $\Gamma =
\Delta_{M_{\alpha_1}}\djoin\ldots\djoin\Delta_{M_{\alpha_r}}$,§ if
$N=\{x^{\alpha_1},\ldots,x^{\alpha_r}\}$. The point here is that $M_N
= M_{\alpha_1}\cup\ldots\cup M_{\alpha_r}$ is the decomposition of
$M_N$ into its connected components: every $M_{\alpha_i}$ is connected
because $x^{\alpha_i}\in\cL$, and since $N$ is discrete, there are no
edges between $M_{\alpha_i}$ and $M_{\alpha_j}$ if $i\ne j$. Therefore
$$\Delta_{M_N}' = \Delta_{M_{\alpha_1}} \djoin \ldots \djoin
\Delta_{M_{\alpha_r}} = \Gamma.$$

For any $N\in D(\cL)$, the set $M_N$ is obviously saturated in
$M$. Conversely, for any saturated subset $S$ of $M$, let $S =
S_1\cup\ldots\cup S_r$ be the decomposition of $S$ into connected
components. Then $N = \{m_{S_1},\ldots,m_{S_r}\} \in D(L_I)$ and since
$S_i = M_{m_{S_i}}$, as $S$ is saturated, it follows that
$M_{m_{S_i}}$ is connected for each $i$, so that $N\in D(\cL)$. This
sets up a one-to-one correspondence between $\sat{M}$ and
$D(\cL)$. Furthermore, under this correspondence $m_S = m_N$ and
$\components{S} = |N|$, so it translates our formula to what we want:
$$b_R =
1+\sum_{S\in\sat{M}}m_S(-z)^{\components{S}+2}\rH(\Delta_S';k)(z).$$
\end{proof}

Fix the field $k$ and introduce the auxiliary notation
$$F(M) = 1+\sum_{S\in \sat{M}} m_S (-z)^{\components{S}+2}
\rH(\Delta_S';k)(z),$$ when $M$ is a set of monomials of degree at
least $2$. If $I$ is a monomial ideal in some polynomial ring $Q$ over
$k$, then set $F(I) = F(M_I)$. So far we have proved that $b_{Q/I} =
F(I)$ whenever $I$ is generated by square-free monomials. The claim of
Theorem \ref{hauptsatz} is that $b_{Q/I} = F(I)$ for all monomial
ideals $I$.

\begin{lemma} \label{F-lemma}
Let $I$ and $I'$ be equivalent monomial ideals, and let $f\colon
L_I\rightarrow L_{I'}$ be an equivalence. Then
$$f(F(I)) = F(I'),$$ where $f(F(I))$ denotes the result of applying
$f$ to the coefficients $m_S$ of $F(I)$, regarding it as a polynomial
in $z$.
\end{lemma}

\begin{proof}
By Lemma \ref{equivalent-antichains}, $f$ induces a bijection of
$\sat{M_I}$ onto $\sat{M_{I'}}$, mapping $S$ to $f(S)$, such that
$\Delta_S'\cong\Delta_{f(S)}'$ and $\components{S} =
\components{f(S)}$ for $S\in\sat{M}$. Since $f(m_S) = m_{f(S)}$ for
all $S\subseteq M_I$, the result follows.
\end{proof}

\begin{proof}[Proof of Theorem \ref{hauptsatz}. General case]
We invoke the construction of Fr\"oberg, \cite{froberg} pp. $30$, in
order to reduce to the square-free case. Let $I$ be any monomial ideal
in $Q=k[x_1,\ldots,x_t]$, and let $M=M_I$. Let $d_i=\max_{m\in
M}\deg_{x_i}(m)$. To each $m\in M$ we associate a square-free monomial
$m'$ in $Q'= k[x_{i,j}\mid 1\leq i\leq t, 1\leq j\leq d_i]$ as
follows: If $m=x_1^{\alpha_1}\cdot\ldots\cdot x_t^{\alpha_t}$ then
$$m' = \prod_{i=1}^t\prod_{j=1}^{\alpha_i}x_{i,j}.$$ The set $M' =
\set{m'}{m\in M}$ minimally generates an ideal in $Q'$, which we
denote by $I'$. The map $M'\rightarrow M$, $m'\mapsto m$, extends to a
map $f\colon L_{I'}\rightarrow L_I$ characterized by the property that
$x_{i,j}$ divides $m\in L_{I'}$ if and only if $x_i^j$ divides
$f(m)$. From this defining property it is easily seen that $f$ is an
equivalence. Therefore $f(F(I')) = F(I)$, by Lemma \ref{F-lemma}.

Let $R$ and $R'$ be the monomial rings associated to $M$ and $M'$
respectively. Using the technique of \cite{froberg} it is easy to see
that
$$b_R(x_1,\ldots,x_t,z) =
b_{R'}(x_1,\ldots,x_1,x_2,\ldots,x_2,\ldots,x_t,z),$$ that is,
$$b_R = f(b_{R'}).$$ But $I'$ is generated by square-free monomials,
so $b_{R'} = F(I')$, whence
$$b_R = f(b_{R'}) = f(F(I')) = F(I),$$
which proves Theorem \ref{hauptsatz} in general.
\end{proof}

\section{Applications and remarks} \label{applications}
We will here give the proofs of the corollaries to the main theorem
and make some additional remarks.

\begin{cor}
With notations as in Theorem \ref{hauptsatz}
$$\deg b_R(z) \leq n+g,$$ where $b_R(z) = b_R(1,\ldots,1,z)$,
$n=|M_I|$ is the number of minimal generators of $I$ and $g$ is the
independence number of $M_I$, i.e., the largest size of a discrete
subset of $M_I$. In particular
$$\deg b_R(z) \leq 2n,$$ with equality if and only if $R$ is a
complete intersection.
\end{cor}

\begin{proof}
If $\Delta$ is a simplicial complex with $v$ vertices, then $\deg
\rH(\Delta;k)(z)\leq v-2$, because either $\dim \Delta = v-1$, in
which case $\Delta$ is the $(v-1)$-simplex and $\rH(\Delta;k) = 0$, or
else $\dim \Delta \leq v-2$, in which case $\rH_i(\Delta;k) = 0$ for
$i> v-2$. The simplicial complex $\Delta_S'$ has $|S|$ vertices. Thus
the $z$-degree of a general summand in the formula \eqref{hauptformel}
for $b_R(\boldx,z)$ is bounded above by $\components{S}+2+|S|-2\leq
g+n$, because the number of components of $S$ can not exceed the
independence number of $M_I$.

Since $g\leq n$ we get in particular that
$$\deg b_{R}(z) \leq 2n,$$ with equality if and only if $M_I$ is
discrete itself, which happens if and only if $R$ is a complete
intersection.
\end{proof}

Now that we know that $R$ and $R'$ below satisfy $b_R = F(I)$ and
$b_{R'} = F(I')$, the next corollary is merely a restatement of Lemma
\ref{F-lemma}.
\begin{cor}
Let $I$ and $I'$ be ideals generated by monomials of degree at least
$2$ in the rings $k[\boldx]$ and $k[\boldx']$ respectively, where
$\boldx$ and $\boldx'$ are finite sets of variables. Let
$R=k[\boldx]/I$, $R' = k[\boldx']/I'$. If $f:L_I\rightarrow L_{I'}$ is
an equivalence, then
$$b_{R'}(\boldx',z) = f(b_R(\boldx,z)),$$ where $f(b_R(\boldx,z))$
denotes the result of applying $f$ to the coefficients of
$b_R(\boldx,z)$, regarding it as a polynomial in $z$. \hfill $\square$
\end{cor}

\begin{remark}
Given formula \eqref{hauptformel}, it is easy to reproduce the result,
implicit in \cite{froberg2} and explicit in
\cite{charalambousandreeves}, that
\begin{equation*}
b_R(\boldx,z) = \sum_{S\subseteq M_I}
(-1)^{\components{S}}z^{|S|+\components{S}}m_S,
\end{equation*}
when the Taylor complex on $M_I$ is minimal. The Taylor complex is
minimal precisely when $m_T=m_S$ implies $S=T$, for $S,T\subseteq
M_I$, i.e., when $L_I$ is isomorphic to the boolean lattice of subsets
of $M_I$. In this case every non-empty subset $S$ of $M_I$ is
saturated, because $m\mid m_T$ implies $m\in T$ for any $T\subseteq
M_I$, and $\Delta_S'$ is a triangulation of the $(|S|-2)$-sphere,
because $m_S=m_M$ only if $S=M$.
\end{remark}

\begin{remark}
The set of saturated subsets of $M$ constitutes a lattice with
intersection as meet and the saturation of the union as join. For each
$m\in M$, the singleton $\{m\}$ is saturated, and the set
$M'=\set{\{m\}}{m\in M}$ is a \emph{cross-cut} of the lattice
$\satt{M}$, in the sense of \cite{folkman}, that is, it is a maximal
antichain. A subset $S\subseteq M'$ is said to \emph{span} if its
supremum in $\satt{M}$ is $M$ and if its infimum is $\emptyset$. The
simplicial complex of subsets of $M'$ that do not span can be
identified with the complex
$$\Gamma_M = \set{S\subseteq M}{\saturation{S}\ne M}.$$ Therefore, by
\cite{folkman} Theorem $3.1$, there is an isomorphism
$$\H(\Gamma_M) \cong \H(\sattt{M}),$$ where $\sattt{M}$ is the
partially ordered set $\satt{M}\setminus\{\emptyset,M\}$. As usual,
the homology of a partially ordered set $P$ is defined to be the
homology of the simplicial complex of chains in $P$.

Actually, one can check that $\saturation{S}\ne M$ is equivalent to
$m_S\ne m_M$ or $S\cap M_i$ disconnected for some connected component
$M_i$ of $M$, so in fact we have an equality $\Gamma_M =
\Delta_M'$. Thus $\Delta_M'$ computes the homology of
$\sattt{M}$. Moreover, if $S\in\satt{M}$ then $\satt{S}$ is equal to
the sublattice $\satt{M}_{\subseteq S} = \set{T\in\satt{M}}{T\subseteq
S}$ of $\satt{M}$. Therefore the homology groups occuring in
\eqref{hauptformel} can be interpreted as the homology groups of the
lower open intervals $\sattt{S} = (\emptyset,S)_{\satt{M}}$ of the
lattice $\satt{M}$ and we may rewrite \eqref{hauptformel} as
\begin{equation} \label{hauptformel-lattice}
b_R(\boldx,z) = 1+\sum_{S\in\sat{M}}
m_S(-z)^{\components{S}+2}\rH((\emptyset,S)_{\satt{M}};k)(z).
\end{equation}
This could be compared with the result of \cite{gasharov} that the
Betti numbers of a monomial ring $R = k[x_1,\ldots,x_t]/I$ can be
computed from the homology of the lower intervals of the lcm-lattice,
$L_I$, of $I$. Specifically, Theorem $2.1$ of \cite{gasharov} can be
stated as
\begin{equation} \label{betti}
\P_R^Q(\boldx,z) = 1 + \sum_{1\ne m\in L_I} m z^2
\rH((1,m)_{L_I};k)(z).
\end{equation}
Here $\P_R^Q(\boldx,z)$ is the polynomial
$$\P_R^Q(\boldx,z) = \sum_{i\geq 0, \alpha\in\zn} \dim_k
\Tor_{i,\alpha}^Q(R,k) x^\alpha z^i.$$

The lattices $\satt{M}$ and $L_M$ are not unrelated. There is a
surjective morphism of join-semilattices from $\satt{M}$ to $L_M$
sending $S$ to $m_S$, so $L_M$ is always a quotient semilattice of
$\satt{M}$. If the graph of $M$ is complete, then the morphism is an
isomorphism.
\end{remark}

As a conclusion, we note how our formula gives a combinatorial
criterion for when a monomial ring is Golod. Interesting sufficient
combinatorial conditions have been found earlier, see for instance
\cite{herzog}, but the author is not aware of any necessary condition
which is formulated in terms of the combinatorics of the monomial
generators.

Recall that $R$ is called a \emph{Golod ring} if there is an equality
of formal power series
$$\P_k^R(\boldx,z) =
\frac{\prod_{i=1}^t(1+x_iz)}{1-z(\P_R^Q(\boldx,z)-1)}.$$ In terms of
the denominator polynomial the condition reads
\begin{equation} \label{golod eq}
b_R(\boldx,z) = 1-z(\P_R^Q(\boldx,z)-1).
\end{equation}

It is easily seen that $S$ is saturated in $M$ if and only if $S$ is
saturated in $M_{m_S}$. Note also that $(1,m)_{L_M} =
L_{M_m}\setminus\{1,m\} =: \bar{L}_{M_m}$. Therefore, after plugging
the formulas \eqref{hauptformel-lattice} and \eqref{betti} into
\eqref{golod eq} and equating the coefficients of each $m\in L_I$, we
get a criterion for $R$ to be a Golod ring as follows:

Call a monomial set $N$ \emph{pre-Golod over $k$} if
$$\rH(\bar{L}_N,k)(z) = \sum_{\substack{S\in\sat{N} \\
m_S=m_N}}(-z)^{\components{S}-1}\rH((\emptyset,S)_{\satt{N}};k)(z).$$
\begin{theorem}
Let $k$ be a field and let $I$ be a monomial ideal in
$k[x_1,\ldots,x_t]$ with minimal set of generators $M$. Then the
monomial ring $R = k[x_1,\ldots,x_t]/I$ is Golod if and only if every
non-empty subset of $M$ of the form $M_m$, with $m§\in L_I$, is
pre-Golod over $k$. \hfill $\square$
\end{theorem}

\bibliographystyle{amsplain}
\bibliography{library}

\providecommand{\bysame}{\leavevmode\hbox to3em{\hrulefill}\thinspace}
\providecommand{\MR}{\relax\ifhmode\unskip\space\fi MR }
% \MRhref is called by the amsart/book/proc definition of \MR.
\providecommand{\MRhref}[2]{%
  \href{http://www.ams.org/mathscinet-getitem?mr=#1}{#2}
}
\providecommand{\href}[2]{#2}
\begin{thebibliography}{10}

\bibitem{avramov-inf}
L.~L. Avramov, \emph{Infinite free resolutions}, Six lectures on commutative
  algebra (Bellaterra, 1996), Progr. Math., vol. 166, Birkh\"auser, 1998,
  pp.~1--118.

\bibitem{avramov}
\bysame, \emph{Homotopy {L}ie algebras and {P}oincar\'e series of algebras with
  monomial relations}, Homology Homotopy Appl. \textbf{4} (2002), no.~2,
  17--27, The Roos Festschrift, vol. 1.

\bibitem{backelin}
J.~Backelin, \emph{Les anneaux locaux \`a relations monomiales ont des s\'eries
  de {P}oincar\'e-{B}etti rationelles}, C.R. Acad. Sc. Paris (1982), 607--610.

\bibitem{bruns}
W.~Bruns and J.~Herzog, \emph{{C}ohen-{M}acaulay rings}, second ed., Cambridge
  Studies in Advanced Mathematics, no.~39, Cambridge University Press, 1996.

\bibitem{charalambousandreeves}
H.~Charalambous and A.~Reeves, \emph{{P}oincar\'e series and resolutions of the
  residue field over monomial rings}, Comm. Algebra \textbf{23} (1995),
  2389--2399.

\bibitem{halperin}
Y.~F\'elix, S.~Halperin, and J-C. Thomas, \emph{Rational homotopy theory},
  Grad. Texts in Math., no. 205, Springer Verlag, 2000.

\bibitem{folkman}
J.~Folkman, \emph{The homology groups of a lattice}, J. Math. Mech. \textbf{15}
  (1966), 631--636.

\bibitem{froberg2}
R.~Fr\"oberg, \emph{Some complex constructions with applications to
  {P}oincar\'e series}, Lecture Notes in Math. \textbf{740} (1979), 272--284.

\bibitem{froberg}
\bysame, \emph{A study of graded extremal rings and of monomial rings}, Math.
  Scand. \textbf{51} (1982), 22--34.

\bibitem{gasharov}
V.~Gasharov, I.~Peeva, and V.~Welker, \emph{The {LCM}-lattice in monomial
  resolutions}, Math. Res. Lett. \textbf{6} (1999), 521--532.

\bibitem{gulliksen-levin}
T.H. Gulliksen and G.~Levin, \emph{Homology of local rings}, no.~20, Queen's
  Papers Pure Appl. Math., 1969.

\bibitem{herzog}
J.~Herzog, V.~Reiner, and V.~Welker, \emph{Componentwise linear ideals and
  {G}olod rings}, Michigan Math. J. \textbf{46} (1999), no.~2, 211--223.

\end{thebibliography}
\end{document}